\documentclass[letterpaper]{amsart}

\newcommand{\myauthor}{Benjamin Antieau and Ben Williams}
\newcommand{\mytitle}{On the classification of oriented $3$-plane bundles over a $6$-complex}

\title[$\SO_3$-bundles]{On the classification of oriented $3$-plane bundles over a $6$-complex}
\author{\myauthor}

\usepackage[pdfstartview=FitH,
            pdfauthor={\myauthor},
            pdftitle={\mytitle},
            colorlinks,
            linkcolor=reference,
            citecolor=citation,
            urlcolor=e-mail,
            ]{hyperref}

\usepackage{amsmath}
\usepackage{amscd}
\usepackage{amsbsy}
\usepackage{amssymb}
\usepackage{verbatim}
\usepackage{eufrak}
\usepackage{microtype}
\usepackage{hyperref}
\usepackage{mathrsfs}
\usepackage{amsthm}
\usepackage[all,cmtip]{xy}
\usepackage{tikz}
\usetikzlibrary{matrix,arrows}
\usepackage[abbrev,lite]{amsrefs}

\usepackage{mathpazo}

\usepackage{color}
\definecolor{todo}{rgb}{1,0,0}
\definecolor{conditional}{rgb}{0,1,0}
\definecolor{e-mail}{rgb}{0,.40,.80}
\definecolor{reference}{rgb}{.20,.60,.22}
\definecolor{mrnumber}{rgb}{.80,.40,0}
\definecolor{citation}{rgb}{0,.40,.80}

\DeclareMathAlphabet{\mathpzc}{OT1}{pzc}{m}{it}
\usepackage{dsfont}

\DeclareMathOperator{\BSU}{BSU}

\DeclareMathOperator{\PU}{PU}

\DeclareMathOperator{\GL}{GL}

\DeclareMathOperator{\Hoh}{H}

\DeclareMathOperator{\sk}{sk}
\DeclareMathOperator{\Sq}{Sq}

\DeclareFontEncoding{OT2}{}{}

\newcommand{\we}{\simeq}
\newcommand{\iso}{\cong}

\newcommand{\ZZ}{\mathds{Z}}

\let\oldmarginpar\marginpar
\renewcommand\marginpar[1]{\-\oldmarginpar[\raggedleft\footnotesize #1]%
{\raggedright\footnotesize #1}}

\newcommand{\SU}{\mathrm{SU}}
\newcommand{\SO}{\mathrm{SO}}
\newcommand{\BSO}{\mathrm{BSO}}

\theoremstyle{plain}
\newtheorem{theorem}{Theorem}

\newtheorem{corollary}[theorem]{Corollary}

\newtheoremstyle{named}{}{}{\itshape}{}{\bfseries}{.}{.5em}{#1 \thmnote{#3}}
\theoremstyle{named}

\theoremstyle{definition}
\newtheorem{definition}[theorem]{Definition}

\theoremstyle{remark}

\begin{document}

\maketitle

\begin{abstract}
    In this short note, we complete the description of low-degree characteristic classes
    of oriented $3$-plane bundles over a $6$-complex. Our goal is to 
    point out and correct an error in L. M. Woodward's 1982 paper ``The classification of
    principal $\PU_n$-bundles over a $4$-complex.''

    \paragraph{\bf Key Words}
    Classification of oriented vector bundles, Postnikov towers.

    \paragraph{\bf Mathematics Subject Classification 2010}
    Primary:
    \href{http://www.ams.org/mathscinet/msc/msc2010.html?t=55Rxx&btn=Current}{55R10},
    \href{http://www.ams.org/mathscinet/msc/msc2010.html?t=55Rxx&btn=Current}{55R45}.
\end{abstract}

When $G$ is a compact Lie group,
the problem of classifying principal $G$-bundles over a topological space $X$ is one of the
central problems of algebraic topology. Attempts to solve this problem in various cases
lead to the ideas of
Postnikov towers, characteristic classes, and $K$-theory. In algebraic geometry, the
analogous problem of understanding vector bundles (principal $\GL_n$-bundles) is the heart
of an enormous amount of ongoing conjecture and research.

When $G=\SO_n$, the problem is to classify oriented $n$-plane bundles on a $k$-complex $X$.
A typical approach is to consider $k$-tuples of characteristic classes of oriented $n$-plane bundles
and to ask which $k$-tuples of cohomology classes of $X$ occur as the characteristic classes
of an oriented $n$-plane bundle. This subject has been
studied since the beginnings of algebraic topology; see~\cite{cadek-vanzura} for an
overview and~\cite{dold-whitney} for early work on this problem. Most results are for small $k$, when $n$ is large with respect to $k$. This is
thanks to simplifications due to Bott periodicity (see~\cite{woodward-orientable}). Another case that has been studied for
$k\leq 9$ is when $k=n$. The case when $n<k$ is much, much more difficult. When $n=2$, the
answer is known: an oriented $2$-plane bundle is determined by its Euler class in
$\Hoh^2(X,\ZZ)$, and any class in $\Hoh^2(X,\ZZ)$ is the Euler class of a $2$-plane
bundle. When $n=3$, previous results allow for a classification when $k\leq 5$
by~\cite{cadek-vanzura}*{Theorem 2}.

\begin{theorem}
    Let $X$ be a $6$-dimensional CW complex. Consider the map
    \begin{equation}\label{eq:5}
        [X,\BSO_3]\rightarrow\Hoh^2(X,\ZZ/2)\times\Hoh^4(X,\ZZ)
    \end{equation}
    which sends an oriented
    $3$-plane bundle $\xi$ on $X$ to the pair of cohomology classes $(w_2(\xi),p_1(\xi))$, where
    $w_2(\xi)$ is the second Stiefel-Whitney class and $p_1(\xi)$ is the first Pontrjagin
    class. Let $\rho_{4*}:\Hoh^4(X,\ZZ)\rightarrow\Hoh(X,\ZZ/4)$ denote reduction modulo
    $4$, and let $P_2:\Hoh^2(X,\ZZ/2)\rightarrow\Hoh^4(X,\ZZ/4)$ be the Pontrjagin square.
    The image of~\eqref{eq:5} consists of the set of classes $(x,y)$
    satisfying
    \begin{equation}\label{eq:1}
        \rho_{4*}(p_1(\xi))=3P_2(w_2(\xi)),
    \end{equation}
    such that $u(x,y)=0$, where $u(x,y)$ is a certain function on the set of pairs of cohomology
    classes satisfying~\eqref{eq:1} with values in $\Hoh^6(X,\ZZ/2)$, to
    be defined below. Moreover, there is a $6$-dimensional CW complex $X$ and a pair of classes
    $(x,y)$ satisfying~\eqref{eq:1} such that $u(x,y)\neq 0$.
\end{theorem}

Previous work of Woodward~\cite{woodward},
using the language of principal $\PU_2$-bundles and the exceptional
isomorphism $\PU_2\iso\SO_3$\footnote{In fact, it was the classification of principal
$\PU_2$-bundles which originally sparked our interest in this problem.}, purported to solve
this problem as well. Unfortunately, the part of the main theorem of Woodward
dealing with $3$-planes on $6$-complexes is mistaken, because of the incorrect assumption
there that $\pi_5\BSO_3=0$, which appears on p.521. As shown by Bott~\cite{bott}*{Theorem 5},
$\pi_5\BSO_3=\ZZ/2$. We explain how this affects the main theorem of~\cite{woodward}, and how to
correct the theorem. The correction requires more than simply re-writing Woodward's proof to
take the correct homotopy group into account; we need additional information, which comes
from the Postnikov tower of $\BSO_3$.

Woodward claimed that the image of the map~\eqref{eq:5} consists of classes
satisfying~\eqref{eq:1} with no other restrictions when $X$ is a CW complex of dimension at
most $6$.
Only a small portion of this claim is false: when $\dim X=6$, there are some classes
$(x,y)\in\Hoh^2(X,\ZZ/2)\times\Hoh^4(X,\ZZ)$ satisfying $\rho_{4*}(y)=3P_2(x)$, but which are
not the characteristic classes of any $\PU_2$-bundle over $X$. The necessary additional
condition is as stated in our theorem.


For $6$-complexes $X$, there is a surjection $[X,\BSO_3]\rightarrow [X,\BSO_3[5]]$, where
$\BSO_3[5]$ denotes the $5$th stage in the Postnikov tower for $\BSO_3$.
The characteristic classes above are obtained by
showing that $\BSO_3[4]$ is equivalent to the homotopy fiber of the map $-3P_2+\rho_{4*}$
\begin{equation}\label{eq:2}
    K(\ZZ/2,2)\times K(\ZZ,4)\rightarrow K(\ZZ/4,4).
\end{equation}
Thus, given a $3$-plane bundle $\xi$ over $X$, the characteristic classes are given by the
composition $$X\xrightarrow{\xi}\BSO_3\rightarrow\BSO_3[4]\rightarrow K(\ZZ/2,2)\times
K(\ZZ,4).$$ The relation~\eqref{eq:1} is expressed in the fact that this map lands in the
fiber of~\eqref{eq:2}. We are left with the problem of computing the image of
$[X,\BSO_2[5]]\rightarrow [X,\BSO_2[4]]$. 

The $5$th stage of the Postnikov tower for $\BSO_3$ gives a fiber sequence,
\begin{equation*}
    K(\ZZ/2,5)\rightarrow\BSO_3[5]\rightarrow\BSO_3[4],
\end{equation*}
which deloops to a map $\BSO_3[4]\rightarrow K(\ZZ/2,6)$ of which $\BSO_3[5]$ is the
homotopy fiber. Under the correspondence between maps $\BSO_3[4]\rightarrow K(\ZZ/2,6)$
and cohomology classes in $\Hoh^6(\BSO_3[4],\ZZ/2)$, the map
is classified by a class $u\in\Hoh^6(\BSO_3[4],\ZZ/2)$. 

\begin{definition}
    Given a space $X$ and classes
    $(x,y)\in\Hoh^2(X,\ZZ/2)\times\Hoh^4(X,\ZZ)$ satisfying~\eqref{eq:1}, one has a uniquely
    determined map $f:X\rightarrow\BSO_3[4]$, and hence a cohomology class
    $f^*(u)\in\Hoh^6(X,\ZZ/2)$. In order
    for $f$ to lift to a map $X\rightarrow\BSO_3[5]$, it is necessary for $f^*(u)=0$. If $\dim
    X\leq 6$, this is also a sufficient condition. If $f$ is determined by classes $(x,y)$
    satisfying~\eqref{eq:1},
    write $u(x,y)$ for $f^*(u)$; thus, $u(x,y)\in\Hoh^6(X,\ZZ/2)$.
\end{definition}

\begin{proof}[Proof of theorem]
    Since Woodward identified $\BSO_3[4]$ as the fiber of the map $K(\ZZ/2,2)\times
    K(\ZZ,4)\rightarrow K(\ZZ/4,4)$, given by relation~\eqref{eq:1}, the image of the
    composition
    \begin{equation*}
        [X,\BSO_3]\rightarrow[X,\BSO_3[4]]\rightarrow\Hoh^2(X,\ZZ/2)\times \Hoh^4(X,\ZZ)
    \end{equation*}
    consists of pairs of classes satisfying~\eqref{eq:1}. By the theory of Postnikov towers,
    a map $f:X\rightarrow\BSO_3[4]$ lifts to $X\rightarrow\BSO_3[5]$ if and only if
    $f^*(u)=0$, which by our definition, occurs if and only if $u(x,y)=0$. Since on a
    $6$-complex any map $X\rightarrow\BSO_3[5]$ lifts to a map $X\rightarrow\BSO_3$, this proves the
    first statement.

    To prove the second statement is
    equivalent to showing that the extension $K(\ZZ/2,5)\rightarrow\BSO_3[5]\rightarrow
    \BSO_3[4]$ is non-split. Indeed, if it is non-split, then the $6$-skeleton of
    $\BSO_3[4]$ together with the composition
    $$\sk_6\left(\BSO_3[4]\right)\rightarrow\BSO_3[4]\rightarrow K(\ZZ/2,2)\times
    K(\ZZ,4)$$ gives an example.

    Recall that $\PU_2=\SU_2/\ZZ/2$, where $\ZZ/2$ is the center of the special unitary
    group $\SU_2$.
    The quotient map $\SU_2\rightarrow\PU_2\iso\SO_3$ induces a map on classifying spaces
    $\BSU_2\rightarrow\BSO_3$, which induces an isomorphism on homotopy
    groups $\pi_i$ for $i>2$. By the naturality of Postnikov towers, there is thus a
    map of extensions
    \begin{equation*}
        \begin{CD}
            K(\ZZ/2,5)  @>>>    \BSU_2[5]   @>>>    K(\ZZ,4)\\
            @|                  @VVV                @VVV\\
            K(\ZZ/2,5)  @>>>    \BSO_3[5]   @>>>    \BSO_3[4].
        \end{CD}
    \end{equation*}
    If the class of the extension in $\Hoh^6(K(\ZZ,4),\ZZ/2)$ is non-zero, then by the
    commutativity of the diagram, the class in $\Hoh^6(\BSO_3[4],\ZZ)$ is non-zero. It
    is not hard to show, using the Serre spectral sequence, that
    $\Hoh^6(K(\ZZ,4),\ZZ/2)=\ZZ/2$, generated by a class $\gamma$. On the other hand, $\Hoh^*(\BSU_2,\ZZ)=\ZZ[c_2]$,
    where the class $c_2$ has degree $4$. Therefore, $\Hoh^6(\BSU_2,\ZZ/2)=0$. Since
    $\BSU_2\rightarrow\BSU_2[5]$ is a $6$-equivalence, it follows that
    $\Hoh^6(\BSU_2[5],\ZZ/2)=0$ as well. If the extension were split, then the
    pullback of $\gamma$ to $\BSU_2[5]$ would be non-zero. Thus the extension is not
    split.
\end{proof}

In~\cite{aw4}, we produce an example of a $6$-dimensional smooth affine variety $X$ and a fixed
non-zero class $x\in\Hoh^2(X,\ZZ/2)$ such that there is no oriented $3$-plan $\xi$ with
$w_2(\xi)=x$. This is despite the fact that there is a pair $(x,y)$ satisfying~\eqref{eq:1}.
Thus, in some sense, Woodward's statement can fail as badly as possible in some situations.


Now, we prove a corollary, which amounts to determining the class $u$ in
$$\Hoh^6(\BSO_3[4],\ZZ/2).$$ By Serre~\cite{serre}*{Section 9}, the $\ZZ/2$-cohomology of
$K(\ZZ/2,2)$ is a polynomial ring
\begin{equation*}
    \Hoh^*(K(\ZZ/2,2),\ZZ/2)=\ZZ/2[u_2,\Sq^1u_2,\Sq^2\Sq^1u_2,\ldots,\Sq^{2^k}\Sq^{2^{k-1}}\cdots\Sq^2\Sq^1u_2,\ldots],
\end{equation*}
where $u_2$ is the fundamental class in degree $2$, and $\Sq^i$ denotes the $i$th Steenrod
operation. Let $\BSO_3[4]\rightarrow K(\ZZ/2,2)$ be denoted by $p$.

\begin{corollary}
    The set $\{u,p^*u_2^3,p^*(\Sq^1u_2)^2\}$ forms a basis of the $3$-dimensional
    vector space $\Hoh^6(\BSO_3[4],\ZZ/2)$.
    \begin{proof}
        There is an isomorphism $\Hoh^*(\BSO_3,\ZZ/2)\iso\ZZ/2[w_2,w_3]$, where $w_i$ has degree $i$. The map
        $\BSO_3\rightarrow\BSO_3[2]\we K(\ZZ/2,2)$ is a $4$-equivalence, so that
        $w_2$ is the pullback of $u_2$, and $w_3$ is the pullback of $\Sq^1u_2$. It follows,
        in fact, that $\Hoh^*(\BSO_3[n],\ZZ/2)$ contains the algebra $\ZZ/2[w_2,w_3]$ for
        $n\geq 2$. A brief examination of the Serre spectral sequence for the fibration
        $K(\ZZ,4)\rightarrow\BSO_3[4]\rightarrow K(\ZZ/2,2)$ shows that the dimension of
        $\Hoh^6(\BSO_3[4],\ZZ/2)$ is at most $3$. The classes $w_2^3$ and
        $w_3^2$ must survive and be distinct, since they do in the cohomology of $\BSO_3$.
        Finally, since we showed in the proof of the theorem that the extension class $u$
        restricts to the non-zero class in $\Hoh^6(K(\ZZ,4),\ZZ/2)$, it follows that the
        asserted classes form a basis for $\Hoh^6(\BSO_3[4],\ZZ/2)$, as desired.
    \end{proof}
\end{corollary}


\begin{bibdiv}
\begin{biblist}


\bib{aw3}{article}{
    author = {Antieau, Benjamin},
    author = {Williams, Ben},
    title = {The topological period-index problem for $6$-complexes},
    journal = {J. Top.},
    eprint = {http://arxiv.org/abs/1208.4430},
    note = {To appear.}
}

\bib{aw4}{article}{
    author = {Antieau, Benjamin},
    author = {Williams, Ben},
    title = {Unramified division algebras do not always contain Azumaya maximal orders},
    journal = {Inv. Math.},
    eprint = {http://dx.doi.org/doi:10.1007/s00222-013-0479-7},
    year = {2013},
}

%




\bib{bott}{article}{
    author={Bott, Raoul},
    title={The space of loops on a Lie group},
    journal={Michigan Math. J.},
    volume={5},
    date={1958},
    pages={35--61},
    issn={0026-2285},
}


\bib{cadek-vanzura}{article}{
    author={\v{C}adek, Martin},
    author={Van{\v{z}}ura, Ji{\v{r}}{\'{\i}}},
    title={On the classification of oriented vector bundles over $5$-complexes},
    journal={Czechoslovak Math. J.},
    volume={43(118)},
    date={1993},
    number={4},
    pages={753--764},
    issn={0011-4642},
    eprint={http://www.math.muni.cz/~cadek/list.html},
}






\bib{dold-whitney}{article}{
    author={Dold, A.},
    author={Whitney, H.},
    title={Classification of oriented sphere bundles over a $4$-complex},
    journal={Ann. Math.},
    date={1959},
    volume={69},
    pages={667--677},
}










%


\bib{serre}{article}{
    author={Serre, Jean-Pierre},
    title={Cohomologie modulo $2$ des complexes d'Eilenberg-MacLane},
    journal={Comment. Math. Helv.},
    volume={27},
    date={1953},
    pages={198--232},
    issn={0010-2571},
}





\bib{woodward}{article}{
    author={Woodward, L. M.},
    title={The classification of principal ${\rm PU}_{n}$-bundles over a $4$-complex},
    journal={J. London Math. Soc. (2)},
    volume={25},
    date={1982},
    number={3},
    pages={513--524},
    issn={0024-6107},
}

\bib{woodward-orientable}{article}{
    author={Woodward, L. M.},
    title={The classification of orientable vector bundles over CW-complexes
    of small dimension},
    journal={Proc. Roy. Soc. Edinburgh Sect. A},
    volume={92},
    date={1982},
    number={3-4},
    pages={175--179},
    issn={0308-2105},
}

\end{biblist}
\end{bibdiv}
\end{document}